\documentclass[a4paper,12pt]{amsart}
\usepackage{amsmath,amssymb}
\newtheorem{theorem}{Theorem}

\newtheorem{lemma}[theorem]{Lemma}
\newtheorem{definition}{Definition}
\newenvironment{proof*}{\vskip 2mm\noindent {}}{\hfill $\Box$ \vskip 2mm}
\def\C{\mathbb C}

\def\P{\mathbb P}
\def\R{\mathbb R}
\def\Z{\mathbb Z}

\def\eps{\varepsilon}

\newcommand{\dist}{\operatorname{dist}}

\title[Projective squeezing function]
{An analogue of the squeezing function for projective maps}

\author{Nikolai Nikolov and Pascal J.~Thomas}

\address{N. Nikolov\\ Institute of Mathematics and Informatics\\Bulgarian Academy
of Sciences\\Acad. G. Bonchev 8, 1113 Sofia, Bulgaria\newline
\indent Faculty of Information Sciences\\
State University of Library Studies and Information Technologies\\
Shipchenski prohod 69A, 1574 Sofia, Bulgaria}\email{nik@math.bas.bg}

\address{P.J. Thomas\\
Institut de Math\'ematiques de Toulouse; UMR5219 \\
Universit\'e de Toulouse; CNRS \\
UPS, F-31062 Toulouse Cedex 9, France} \email{pascal.thomas@math.univ-toulouse.fr}

\keywords{Projective maps, invariant distances, squeezing function}

\subjclass[2010]{52A20, 53A20, 32F45}

\thanks{The first named author is partially supported by the Bulgarian National Science Fund,
Ministry of Education and Science of Bulgaria under contract DN 12/2. This paper was started
while his was visiting the Paul Saba\-tier University, Toulouse in November 2018
as a guest professor.}

\begin{document}

\begin{abstract}
In the spirit of Kobayashi's applications of methods of invariant metrics to questions of projective
geometry, we introduce a projective analogue of the complex squeezing function.
Using Frankel's work, we prove that for convex domains it stays uniformly bounded from below.
In the case of strongly convex domains, we show that it tends to $1$ at the boundary. This is
applied to get a new proof of a projective analogue of the Wong-Rosay theorem.
\end{abstract}

\maketitle

\section{Introduction}

The projective maps are the ones that preserve lines in projective space. They are linear
in the homogeneous coordinates, and in  affine space yield linear-fractional maps (which
we will call projective too, with a slight abuse of language).
There is a long tradition of applying the appropriate analogues of convex objects
to complex analysis.
Surprisingly, it is also sometimes useful to study geometrically
convex domains and projective maps, which are rather rigid objects, with the methods
developed for complex analysis in several variables.
For instance, one can use projective mappings from an interval into domains to
construct metrics analogous to the Carath\'eodory and Kobayashi metrics, which
recover the classical Hilbert metric in the case of convex domains.
Shoshichi Kobayashi developed this approach in \cite{Ko}, and L\'aszl\'o Lempert
summarized the analogy and built upon it in \cite{Le}. This was pursued in papers such as
\cite{Fra} and \cite{Zim}.

The complex squeezing function was defined under this name in \cite{DGZ},
which provides a good overview of the motivations to study it. It has been
the object of numerous further works in recent years.
We will be using a ``projective'' analogue of the squeezing function study its
relationship with the properties of convex sets.
Some of the results of Sidney Frankel's pioneering paper \cite{Fra} can be rephrased
as the fact that the projective squeezing function of convex domains
is bounded from below by a constant depending only on the dimension. We also
give a converse.

We shall use the behavior of the projective squeezing function to give a necessary condition
for a point of the boundary to be strictly convex (see Definition \ref{strcvx}).

We fix some notations.  Let $\P \R^d = \R^{d+1}\setminus \{0\}/\sim$, where
$x\sim y$ means that there exists $\lambda \in \R\setminus \{0\}$ such that $x=\lambda y$.
As usual, embed $\R^d $ as the classes
represented by $\{(1:x_1: \dots :x_n)\}$ in  $\P \R^d $, which form a dense open set.
A \emph{projective map} of $\P \R^d$ is a map induced by a linear map of $\R^{d+1}$.
We assimilate it with its restriction to $\R^d $.

\begin{definition}
Given a domain $D \subset \R^d$ and $z\in D$, the (projective) \emph{squeezing function}
of $D$ at $z$ is
\begin{multline*}
s_D(z):= \sup\left\{ r>0 : \exists \Phi \mbox{ a projective map s. t. } \Phi(z)=0,
\right.
\\
\left.
\Phi (D) \subset B(0,1), \mbox{ and } B(0,r) \subset \Phi(D)
\right\}.
\end{multline*}
\end{definition}

We set $s_D=0$ if the domain is not projectively equivalent to a bounded domain (in which
case the above supremum is over an empty set); at the other extreme, if there is $z\in D$ such that $s_D(z)=1$,
then $D$ is projectively equivalent to the ball. As in the complex case \cite[Theorem 3.1]{DGZ},
one can show that $s_D$ is a continuous function of $z$.

The complex squeezing function was defined by looking at (bounded) domains in $\C^d$
and holomorphic maps. It yields a new holomorphic invariant.

One of the motivations to introduce the squeezing function in the holomorphic case
was to compare the infinitesimal Kobayashi-Royden and Carath\'eodory-Reiffen (pseudo)metrics.
For a point $z$ in a domain $D\subset \C^d$ and a vector $v\in\C^d$, let $\kappa_D(z;v)$ stand
for the Kobayashi-Royden metric and $\gamma_D(z;v)$ stand
for the Carath\'eodory-Reiffen metric, and $S_D(z)$ for the holomorphic squeezing function at $z$,
as defined in \cite{DGZ}. Then one can  show \cite{DGZ}, using the monotonicity properties
of the invariant metrics and their explicit expression in the ball, that
\[
 S_D(z) \gamma_D(z;v) \le S_D(z) \kappa_D(z;v) \le \gamma_D(z;v).
\]
Bounded domains where $s_D \ge c >0$ are called \emph{holomorphic homogeneous regular domains}
\cite{LSY, LSY2}. On those, the Kobayashi, Carath\'eodory and several other invariant metrics are equivalent.
Those domains include several well-known classes: Teichm\"uller spaces,
bounded domains covering compact K\"ahler manifolds, and strictly convex domains with
$\mathcal C^2$-boundary \cite{Yeu}.

S. Kobayashi defines in \cite[(5.4)]{Ko} a projective analogue of the Koba\-yashi-Royden metric.
Let $I$ be the interval $(-1,1)$. If
$p$ is a point of a domain $D \subset \P^d$ (or more generally a manifold $M$ with a projective structure),
$V$ a tangent vector at $p$,
then
\begin{multline*}
F_D(p;X):=
\\
\inf \left\{ 2|V|: f \mbox{ is a projective map }I\rightarrow D, f(0)=p, Df(0)(V)=X \right\}
\end{multline*}
(the factor $2$ is due to normalizations coming from the Poincar\'e or Hilbert metric).
The analogue of the Carath\'eodory-Reiffen metric can likewise be defined as
\begin{multline*}
\label{defcarapro}
C_D(p;X):=
\\
\sup \left\{ 2|Df(p)(X)|: f \mbox{ is a projective map }D\rightarrow I, f(p)=0 \right\}.
\end{multline*}
Then the analogue of the Schwarz lemma \cite[Lemma 2.5]{Ko} implies that
$C_D(p;X)\le F_D(p;X)$ and one has, with the same proof as in the holomorphic case,
\begin{equation}
\label{CKcomp}
s_D (p) F_D(p;X) \le C_D(p;X).
\end{equation}

Let us add that \cite{Ko} also defines the analogues of the Kobayashi and Carath\'eodory
pseudodistances for  projective structures, and proves that for convex domains in $\P\R^d$,
they coincide with the Hilbert pseudodistance \cite[Example 3.17]{Ko}. That last distance is defined
as follows: given two points $p, q \in D$,
\begin{equation}
\label{hilb}
d_D(p,q):= \left| \log (ab;pq) \right|,
\end{equation}
where $a, b$ are the points where the line $\overline{pq}$ crosses $\partial D$
and $(ab;pq) $ denotes the cross ratio of those four points.

Note that $F_D$ is the infinitesimal form of the above pseudodistance:
\begin{equation}
\label{koba}
F_D(p;X)=\frac{1}{P_+}+\frac{1}{P_-},
\end{equation}
where $$P_\pm=\inf\{\lambda>0:p\pm\lambda X\notin D\}.$$
We also point out
that $C_D=C_{\hat D}=F_{\hat D},$
where $\hat D$ is the (open) convex hull of $D.$

As in \cite{Fra}, we say that a convex domain $D \subset \R^d$ (or $\P\R^d$) is \emph{proper}
 if it contains no
affine line.  We state our results for domains in $\R^d$, but they still hold in the projective
space.

A squeezing function varies between $0$ and $1$. For the projective squeezing
function, being bounded away from $0$ is equivalent to proper convexity of the domain.

\begin{theorem}
\label{lowerbound}
\begin{enumerate}
\item
For every $d\in \Z_+^*$, there is $r_d>0$ such that for any proper convex domain $D\subset \R^d$,
for any $z\in D$, $s_D(z) \ge r_d$.
\item
If $D\subset \R^d$ is a domain such that $\inf_{z \in D} s_D(z)>0$, then $D$ is convex
and proper.
\end{enumerate}
\end{theorem}

This result evidences a gap in the behavior of the squeezing function: if it does not tend
to $0$ near the boundary, then it must be bounded below by a universal constant $r_d$, so
$\inf_D s_D$ can never take any value in $(0,r_d)$. It
would be interesting to determine the precise value of $r_d$.

The holomorphic analogue of this result holds for proper convex domains \cite[Theorem 1.1]{Fra}
(see also \cite[Theorem 1.1]{KZ}),
as well as for non-degenerate $\C$-convex domains \cite[Theorem 1]{NA}.

We now turn to cases where the squeezing function approaches $1$, rather than $0$, near some point of $\partial D$.
On a ball (or a domain projectively equivalent to a ball), the squeezing function is obviously identically equal to $1$.
The ball is the simplest example of a strictly convex domain, i.e.
a domain where all boundary points are strictly convex.

\begin{definition}
\label{strcvx}
A point $p\in \partial D$ is called \emph{strictly convex} if $\partial D$
is $\mathcal C^2$-smooth in a neighborhood of $p$, and the restriction to the
tangent hyperplane at $p$ of the Hessian of the defining function is definite positive.
\end{definition}

\begin{theorem}
\label{typetwo}
Let $D\subset \R^n$ be a convex domain,
and $p \in \partial D$. If $p$ is a strictly convex boundary point, then $\lim_{x\to p} s_D(x)=1$.
\end{theorem}

In the holomorphic case,  this is
\cite[Theorem 1.3]{DGZ2} (see also \cite[Theorem 4.1]{KZ}). The proof uses essentially the embedding  in \cite[Theorem 1.1]{DFW}.

The converse is an open and interesting question, to be compared with  results in the holomorphic case. Andrew Zimmer \cite[Theorem 1.7]{Zim2}
proved that for $D \subset \C^d$, $d\ge 2$, a bounded convex domain with $\mathcal C^{2,\alpha}$
boundary, if $s_D$ tends to $1$ at the boundary (and even under a slightly weaker hypothesis),
then $D$ must be strongly pseudoconvex. This was motivated by previous results from other authors
showing that any bounded strongly pseudoconvex domain $D \subset \C^d$, $d\ge 2$ with
$\partial D \in \mathcal C^2$ must have $s_D$ tending to $1$ at the boundary, but that
the converse fails, see the references in \cite{Zim2}.

\section{Proofs}

\subsection{Proof of Theorem 1 (1).}

It follows from:
\begin{theorem}{\cite[Theorem 7.6 \& Corollary 7.8, p. 200]{Fra}}

Let $V:= (-1,\infty)^d \subset \R^d$.

For every $d\in \Z_+^*$, there exists $r'_d >0$ such that for any proper convex domain $D \subset \R^d$ with $0 \in D$,
there is an affine map $A$ such that $B(0,r'_d) \subset D \subset V$ and $A(0)=0$ (i.e. $A$ is in fact linear).
\end{theorem}

Denote $x=(x_1, \dots, x_d)\in \R^d$ and define a projective map from $\P^d$ to itself by its values on
$\R^d\setminus \{\sum_{i=1}^d x_i= -d-1 \}$:
\[
\Phi (x):= \frac1{\sqrt d} \frac{x}{d+1+ \sum_{i=1}^d x_i}.
\]
Clearly $\Phi(0)=0$. We claim that for any $r\in (0,1)$,
\[
B\left(0, \frac{r}{dr+\sqrt d (d+1)}\right) \subset \Phi\left( B(0,r) \right) \subset \Phi\left(V\right) \subset B(0,1),
\]
which will finish the proof, with $r_d= \frac{r'_d}{dr'_d+\sqrt d (d+1)}$.

To get the right-hand inclusion, let $x \in V$, then $d+1+ \sum_{i=1}^d x_i>1$, so if for some $k$,
$x_k <0$, then $0> \frac{x_k}{d+1+ \sum_{i=1}^d x_i} > x_k >-1$.

On the other hand, if $x_k\ge0$, then
\[
0\le \frac{x_k}{d+1+ \sum_{i=1}^d x_i} \le \frac{x_k}{d+1+ x_k + \sum_{i: i\ne k} x_i} \le \frac{x_k}{x_k+2} <1.
\]
In each case, $\left| \Phi (x)_k \right| < 1/\sqrt d$, so $\|\Phi (x)\|<1$.

Conversely, note that
\[
\Phi^{-1}(y)= \frac{\sqrt d (d+1)}{1- \sqrt d \sum_{i=1}^d y_i} y.
\]
Choose $y$ with $(dr+\sqrt d (d+1))\|y\|<r$. Then, since $\sum_{i=1}^d y_i\le \|y\| \sqrt d$,
we have $1- \sqrt d \sum_{i=1}^d y_i > 0$, and
\begin{multline*}
\sqrt d (d+1)\|y\| + r \sqrt d \sum_{i=1}^d y_i < r
\Leftrightarrow
\frac{\sqrt d (d+1)}{1- \sqrt d \sum_{i=1}^d y_i} \|y\| < r.
\end{multline*}
\hfill \qed

\subsection{Proof of Theorem 1 (2).}

First notice that $D$ cannot contain a whole line, because this would make any
projective map $f: D\rightarrow B(0,1)$ degenerate, and so $f(D)$ could not contain a non-trivial ball.

Assume that $D$ is not convex, $d\ge 2$. Then
there is a boundary point $p \in \partial D \cap \hat D$, where we recall that $\hat D$ is the (open) convex hull of $D$;
let $\delta >0$ be such that $B(p,\delta)\subset \hat D$.
Let $D \ni p_k \to p$, and
\newline $X_k := \|p-p_k\|^{-1} (p-p_k)$.
Then
\[
C_D (p_k,X_k)=C_{\hat D } (p_k,X_k)\le C_{B(p,\delta)} (p_k,X_k)\le C<\infty.
\]

To estimate $F_D (p_k,X_k)$ when $t$ is close to $0$, we use the fact that $I$
is projectively isometric to $(0,\infty)$ with $0$ going to $1$. Any projective map $f: (0,\infty) \longrightarrow D$
with $f(1)= p_k$ and $Df(0)(V)=X_k$ must verify $f((0,\infty))\subset p+\R_- X_k$, so
\[
F_D (p_k,X_k) \ge F_{(-\infty,0)}(-\|p-p_k\|,1) = \frac1{\|p-p_k\|}.
\]
Therefore $s_D(p_k) \leq C_D (p_k,X_k)/F_D (p_k,X_k)\le C\|p-p_k\|$ and in particular cannot be bounded below as $p_k \to p$.

{\bf Remarks.}

(1) We have proved a bit more; suppose that $\partial D$ has positive
reach, i.e. that when points are close enough to the boundary, there is a unique
closest point to them on it.  This happens in particular when $\partial D$ is $\mathcal C^{1,1}$-smooth.
Writing $\delta_D(z)$ for the distance of a point
$z\in D$ to $\partial D$, we see that if $\delta_D(z) = o\left( s_D(z) \right)$ (uniformly),
then $D$ must be convex, and by Part (1) of the Theorem,
if it is proper (for example, bounded), then $s_D$ must be in fact bounded below,
by a constant which only depends on the dimension.

(2) Using \cite[Theorem 2.1.27]{Hor}, one could find points tending to the boundary such that the argument
above could be carried out with a constant vector $X$ instead of $X_k$.

\subsection{Proof of Theorem 2.}
Our proof is similar to that carried out for the complex case in \cite{FW},
 where more precise results are obtained in the $\mathcal C^3$- and $\mathcal C^4$-smooth settings.

First pick  coordinates so that $p=0$ and the tangent hyperplane to $\partial D$ at $p$ is
given by $x_1=0$. Locally, $D=\{ (x,x'): x_1+f(x')<0\}$, where $x':=(x_2,\dots,x_d)$, $f(0)=0$,
$Df(0)=0$.
The hypothesis implies that the Hessian of $f$ at $0$ is definite positive, and
since it varies continuously, this is true uniformly for the points of $\partial D$ in a neighborhood $U_1$ of $p$.

Suppose we have a sequence $D \ni q_n \to p$. We may assume that $n$ is large enough so that
the orthogonal projection $p_n$ of $q_n$ to $\partial D$ is well defined and belongs to $U_1$.
Performing a translation and rotation, we may choose new coordinates so that $p_n=0$ and the tangent hyperplane to $\partial D$ at $p_n$ is
given by $x_1=0$. Since $D$ is bounded, we may perform a dilation of the $x_1$ coordinate (with coefficient uniformly bounded above and below w.r.t. $n$)
so that $\inf \{x_1: x \in D\}=-1$. Since the Hessian of $f$ at $0$ is uniformly definite positive, we may
apply a linear map in the $x'$ coordinate in $\R^{d-1}$ which is uniformly bounded above and below so that the Hessian becomes the identity matrix.
Finally, we are reduced to the situation where $q_n= (-\delta, 0,\dots,0)$,  $D=\{ (x,x'): x_1+f(x')<0\}$, $f(x')=\|x'\|^2 + o(\|x'\|^2)$,
and $\delta\asymp \dist(q_n,\partial D)$, with constants uniform in $n$.

For each $\delta >0$, consider the projective map
\[
\phi_\delta : I \times \R^{d-1} \ni y \mapsto \frac1{1+y_1}\begin{pmatrix} \delta (y_1-1) \\ \sqrt \delta y' \end{pmatrix} \in (-\infty,0) \times \R^{d-1}.
\]
It sends $B((0,0); 1)$ to $\{ (x,x'): x_1< -\|x'\|^2\}$ and $(0,0)$ to $(-\delta,0)$. We shall now check that $\phi_\delta^{-1} (D)$ lies
between two balls with radii tending to $1$ as $\delta \to 0$.

\subsubsection*{Estimation from below}

For any $\eps>0$, there exist $ R>0$ (uniform in $n$) such that our domain in the new coordinates contains the lens-shaped set $L_R:= \{ -(1+\eps)R^2 < x_1 < -(1+\eps)\|x'\|^2\}$.
Let us compute $\phi_\delta^{-1}(L_R)$.

Setting $x=\phi_\delta(y)$, first,
\begin{multline*}
x_1 < -(1+\eps)\|x'\|^2 \Leftrightarrow \delta \frac{-1+y_1}{1+y_1}  < -(1+\eps) \delta \frac{\|y'\|^2}{(1+y_1)^2}
\\
\Leftrightarrow
y_1^2 +  (1+\eps)\|y'\|^2 <1.
\end{multline*}
This contains the ball $B(0,(1+\eps)^{-1/2})$. So given any $\eta>0$, we can choose $\eps$ (and therefore $R$) so that
$B(0,1-\eta) \subset B(0,(1+\eps)^{-1/2})$.

Second,
\[
-(1+\eps)R^2 < x_1 \Leftrightarrow -1 < \frac{\delta}{(1+\eps)R^2} \frac{-1+y_1}{1+y_1} \Leftrightarrow y_1 >
-\frac{1-\frac{\delta}{(1+\eps)R^2} }{1+\frac{\delta}{(1+\eps)R^2}}.
\]
This last quantity can be made as close to $-1$ as desired by making $\delta$ small, so once $\eps$ and $R$ are fixed as above,
we can choose $\delta$ so that $-1+\eta > -\frac{1-\frac{\delta}{(1+\eps)R^2} }{1+\frac{\delta}{(1+\eps)R^2}}$,
and $B(0,1-\eta) \subset \phi_\delta^{-1}(L_R) \subset \phi_\delta^{-1} (D)$.

\subsubsection*{Estimation from above}

We first include $D$ into a slightly larger domain with a simpler form.
\begin{lemma}
\label{model}
For any $\eps>0$, there exists $r>0$ such that $f(x')\ge \tilde{f} (x')$, where
\[
\tilde{f} (x'):=
\left\{
\begin{matrix}
(1-\eps) \|x'\|^2, \mbox{ for } \|x'\| \le r
\\
 (1-\eps) r(2\|x'\|-r), \mbox{ for } \|x'\| \ge r.
\end{matrix}
\right.
\]
\end{lemma}
Therefore $D \subset \tilde{D}:=\{x+\tilde{f} (x')<0, -1<x_1\}$. Observe that this is a domain with $\mathcal C^1$-smooth boundary
outside of $\{x_1=-1\}$.

\begin{proof}
By Taylor's formula at order $2$ applied to $f$,
\newline $f(x')=\|x'\|^2(1+\xi_1(x'))$, with $\lim_{x'\to0} \xi_1(x')=0$.

By convexity of $f$, for any $x'$ with $ \|x'\| \ge r$,
\[
f(x') \ge f\left(r\frac{x'}{\|x'\|}\right) + Df\left(r\frac{x'}{\|x'\|}\right) \cdot \left( x'- r\frac{x'}{\|x'\|}\right).
\]
Using the fact that $D^2f(0)=2 Id$ and Taylor's formula at order $1$ applied to $Df$, we have for any $h \in \R^{d-1}$
\[
Df\left(r\frac{x'}{\|x'\|}\right) \cdot h = \left\langle 2r \frac{x'}{\|x'\|} + r \xi_2(r\frac{x'}{\|x'\|}), h\right\rangle,
\]
with $\lim_{r\to0} \| \xi_2(r) \|=0$.  Regrouping the terms, and setting $\tilde{\xi_1}(r):= \min_{\|x'\| \le r} \xi_1(x')$,
we have for $\|x'\| \ge r$
\begin{multline*}
f(x') \ge (1+\tilde{\xi_1}(r))r^2 + 2r \|x'\| - 2 r^2 + r \left\langle r \xi_2(r\frac{x'}{\|x'\|}), \left( x'- r\frac{x'}{\|x'\|}\right)\right\rangle
\\
\ge r (\|x'\|-r) + r \|x'\| + \tilde{\xi_1}(r) r^2 - \tilde \xi_2(r) r \|x'\|,
\end{multline*}
where $ \tilde \xi_2(r):= \max_{x'\neq 0} \xi_2(r\frac{x'}{\|x'\|})$. We may now choose $r>0$ so that
$\tilde{\xi_1}(r)>-\eps/2$ , $\tilde{\xi_2}(r) <\eps/2$ and we obtain the Lemma.
\end{proof}

Since $\overline{\tilde{D}}$ is the convex hull of
\[
\left(\partial \tilde{D} \cap \{x_1=-1\}\right) \cup \left\{ x_1=-(1-\eps) \|x'\|^2,  \|x'\| \le r \right\},
\]
its preimage under $\phi_\delta$ will be obtained by taking the convex hull of the preimages of those two sets. The second
 one is given by
$x_1=-1$ and $(1-\eps) r(2\|x'\|-r) \le 1$, i.e. $ \|x'\| \le \frac12 \left( r+ \frac{1}{(1-\eps) r} \right)$. Computing its preimage,
we find
\[
y_1= \frac{-1+\delta}{1+\delta}, \quad \|y'\| \le \frac{\sqrt \delta}{1+\delta} \left( r+ \frac{1}{(1-\eps) r} \right).
\]
An easy computation shows that then $y_1^2 +  \|y'\|^2 \le 1 + O(\delta)$ (the constants in $O$ depend on $r$).

On the other hand, $\phi_\delta^{-1}\left\{ x_1=-(1-\eps) \|x'\|^2,  \|x'\| \le r \right\}$ is contained in
$y_1^2+(1-\eps) \|y'\|^2=1$, and one can check that locally near the origin $\phi_\delta^{-1} (\tilde D)$ will be inside that ellipsoid.
Notice that when we let $\delta\to0$, the preimage of the neighborhood of the origin in $\partial \tilde{D}$ given by $\|x'\| \le r$
will tend to cover the whole ellipsoid (leaving out the point $(-1,0)$ of course).

So given $\eta>0$, taking first $\eps>0$ small enough, then $\delta>0$ very small, $\phi_\delta^{-1} (D) \subset \phi_\delta^{-1} (\tilde{D}) \subset B(0,1+\eta)$.

\section{An application}

The famous Wong-Rosay theorem (see \cite[Theorem]{Won}, \cite[Proposition]{Ros}) states that if a bounded strictly pseudoconvex domain in $\C^d$ has a
non-compact group of holomorphic automorphism, then $D$ is biholomorphically equivalent to the unit ball. This result has been extended in
\cite[Theorem 3]{Efi} to any domain in $\C^d$ having a strictly pseudoconvex boundary point which is an accumulation point of the group action.

The same is true in $\R^d$ if we replace pseudoconvexity by convexity and holomorphicity by projectivity (see \cite[Th\'eor\`eme 1]{SM},
\cite[Theorem 1.1]{Jo}, \cite[Theorem 3]{Yi}).

Our purpose is to give a short proof of this fact by using the squeezing function.

\begin{theorem}
\label{aut}
Let $D$ be a domain in $\R^d.$ Assume that there exist points $p\in\partial D,$
$q\in D$ and a sequence $(\varphi_j)$ of projective automorphisms of $D$ such that $q_j:=\varphi_j(q)\to p.$
If $p$ is strictly convex, then $D$ is projectively equivalent to the unit ball.
\end{theorem}

The converse is obviously true, since the group of projective automorphisms of the unit ball acts transitively.

\begin{proof} Let $U$ be a neighborhood of $p$ such that $D\cap U$ is convex. Set $D_j=\varphi_j^{-1}(D\cap U).$

\begin{lemma}
\label{exh}
$(D_j)$ is an exhaustion of $D.$
\end{lemma}
Then, similarly to \cite[Theorem 2.1]{DGZ2} in the complex case, we have
$$s_{D_j}(q)\to s_D(q).$$

On the other hand, Theorem \ref{typetwo} implies that
$$s_{D_j}(q)=s_{D\cap U}(q_j)\to 1,$$

So, $s_D(q)=1$ and Theorem \ref{aut} follows.
\medskip

\noindent{\it Proof of Lemma \ref{exh}.} Denote by $k_D$ the ``projective'' Kobayashi metric,
i.e. the integrated form of $F_D.$

Let $K\Subset D$ and
$$M:=\sup_K k_D(q,\cdot)=\sup_{\varphi_j(K)}k_D(q_j,\cdot).$$

Since $p$ is strictly convex, it follows by \eqref{koba}
that there exists a neighborhood of $p$, $V\Subset U$, such that
$$2F_D\ge F_{D\cap U}\quad\mbox{on }(D\cap V)\times\R^d.$$
Hence
$$2\inf_{D\cap\partial V}k_D(q_j,\cdot)\ge \inf_{D\cap\partial V}k_{D\cap U}(q_j,\cdot)=:M_j\quad\mbox{if }q_j\in V.$$

Since $q_i\to p,$ using \eqref{hilb}, we may find a $j_0\in\Bbb N$ such that $M_j>2M$ for $j\ge j_0.$ Then $K\subset D_j$
if $j\ge j_0.$
\end{proof}

{\bf Remark.} Similar arguments imply that Lemma \ref{exh} holds in the complex case if $p$ is a local holomorphic peak function.
Then, as above, the holomorphic analogue of Theorem \ref{typetwo}, namely \cite[Theorem 1.3]{DGZ2}, leads to the already mentioned
complex version of Theorem \ref{aut} (see \cite[Theorem 3]{Efi}).

{}

\end{document}